\documentclass[10pt]{amsart}
\usepackage{epsfig, graphics, psfrag, color}
\usepackage{amscd} 
\renewcommand{\setminus}{{\smallsetminus}}

\newcommand{\bdy}{{\partial}} 

\newcommand{\GA}{{\mathbb{G}_A}}
\newcommand{\GB}{{\mathbb{G}_B}}
\newcommand{\GRA}{{\mathbb{G}'_A}}
\newcommand{\GRB}{{\mathbb{G}'_B}}

\newcommand{\twist}{{\mathrm{tw }(D)}}

\DeclareMathSymbol{\minus} {\mathord}{operators}{"2D} 

\newcommand{\area}{{\rm area}}

\renewcommand{\phi}{{\varphi}}

\newcommand{\HH}{{\mathbb{H}}}
\newcommand{\RR}{{\mathbb{R}}}

\def\co{\colon\thinspace}

\usepackage{amsthm, amsmath}
\usepackage{amssymb} 
\usepackage{amscd}

\usepackage[curve,matrix,arrow,frame,tips]{xy}
\usepackage{verbatim}

\newcommand{\lmin}{{\ell_{\rm min}}}

\newcommand{\ltop}{{\ell_{\rm ext}}}
\newcommand{\lbigon}{{\ell_{\rm in}}}
\DeclareMathSymbol{\minus} {\mathord}{operators}{"2D} 
\newcommand{\vol}{{\rm vol}}
\newcommand{\guts}{{\rm guts}}
\newcommand{\abs}[1]{{\left\vert #1 \right\vert}}
\newcommand{\tw}{{\mathrm{tw }}}

\theoremstyle{plain}
\newtheorem{theorem}{Theorem}[section]
\newtheorem{corollary}[theorem]{Corollary}
\newtheorem{lemma}[theorem]{Lemma}
\newtheorem{prop}[theorem]{Proposition}

\theoremstyle{definition}
\newtheorem{define}[theorem]{Definition}
\newtheorem*{remark}{Remark}
\def\arrowdown#1#2{\Big\downarrow \rlap{$\vcenter{\hbox{$\scriptstyle#2$}}$}
{\hbox to -10pt{\hss{$\vcenter{\hbox{$\scriptstyle#1$}}$}}}}

\begin{document}


\title[Symmetric links and {C}onway sums]{Symmetric links and {C}onway
sums:\\ volume and {J}ones polynomial}

\author[D. Futer]{David Futer}
\author[E. Kalfagianni]{Efstratia Kalfagianni}
\author[J. Purcell]{Jessica S. Purcell}

\address[]{Department of Mathematics, Temple
University,
Philadelphia, PA 19122}

\email[]{dfuter@math.temple.edu}

 \address[]{Department of Mathematics, Michigan State
 University,
 East Lansing, MI, 48824}

 \email[]{kalfagia@math.msu.edu}

\address[]{
Department of Mathematics,
Brigham Young University,
Provo, UT 84602}
\email[]{jpurcell@math.byu.edu }

\thanks{{Futer is supported in part by NSF--RTG grant DMS--0353717.
Kalfagianni is supported in part by NSF--FRG grant
DMS-0456155 and by NSF grant DMS--0805942. Purcell is supported in part by NSF grant DMS-0704359.}}

\thanks{ \today}

\begin{abstract}
We obtain bounds on hyperbolic volume for periodic links and Conway
sums of alternating tangles. For links
that are Conway sums we also bound the hyperbolic volume in terms of
the coefficients of the Jones polynomial.
\end{abstract}

\maketitle


\section{Introduction}\label{sec:intro}
Given a combinatorial diagram of a knot in the 3--sphere, there is an
associated 3--manifold, the \emph{knot complement}, which decomposes
into geometric pieces by work of Thurston \cite{thurston:bulletin}.
A central goal of modern knot theory is to relate this geometric
structure to simple topological properties of the knot and to
combinatorial knot invariants.  To date, there are only a handful of
results along these lines.  Lackenby found bounds on the volume of
alternating links based on the number of twist regions in the link
diagram \cite{lackenby:alt-volume}.  We extended these results to all
links with at least seven crossings per twist region in \cite{fkp-07},
and in \cite{fkp-08} we obtain similar results for links that are
closed 3--braids.  Our method is to apply a result bounding the
change of volume under Dehn filling based on the length of the
shortest filling slope.  In all these cases the relation between twist
number and volume was also important in establishing a \emph{coarse
volume conjecture}: a linear correlation between the coefficients of
the classical Jones polynomial and the volume of hyperbolic links.

In the present paper, we build upon the methods of \cite{fkp-07} as
well as very recent work of Gabai, Meyerhoff, and Milley
\cite{gmm:smallest-cusped}; Agol, Storm, and Thurston
\cite{agol-storm-thurston}; and Agol \cite{agol:two-cusped}.  We use
this work to give explicit estimates on the volume for links with
symmetries of order at least six, and to give estimates on the volume
and coefficients of the Jones polynomial under Conway summation of
tangles.  As in the results above, we obtain explicit linear bounds on
volume in terms of the twist number of a diagram.

\subsection{Links with high order of symmetry}\label{symmetry}
A link $K$ is called \emph{periodic} with period an integer $p>1$ if
there exists an orientation--preserving diffeomorphism $h\co {S}^3
\rightarrow {S}^3$ of order $p$, such that $h(K)=K$ and either  $h$ has fixed points or $h^i$ has no fixed points for all
$0<i<p$.   The solutions to
the Smith conjecture \cite{smith} and the spherical spaceforms
conjecture \cite{morgan-tian} imply that $h$ is conjugate to an element
of $SO(4)$. Thus, if $h$ has no
fixed points, the group generated by $h$ acts freely on $S^3$ and the quotient of $S^3$ is a lens space $L(p,q)$.
Furthermore,  the
quotient of $S^3 \setminus K$ is a link complement in $L(p,q)$.
Otherwise, the orthogonal action conjugate to $h$ must be a $2\pi/p$ rotation
about a great circle $C_h \subset S^3$, and the quotient is still
$S^3$.  When the axis $C_h$ is either a component of $K$ or disjoint
from $K$ (in particular, when $p>2$), the quotient of $K$ is a link
$K' \subset S^3$.

\begin{theorem}\label{thm:periodic-smith}
Let $K$ be a hyperbolic periodic link in $S^3$. Assume that the period
of $K$ is $p\geq 6$, and acts by rotation about an axis $C_h$. Let
$K'$ be the quotient of $K$. Then
$$\vol(S^3\setminus K) \; \geq \: p \, \left(1-\frac{2\sqrt
  2\pi^2}{p^2}\right)^{3/2} 
  \vol(S^3 \setminus K').
$$
\end{theorem}

In the statement above, $S^3\setminus K'$ may or may not be
hyperbolic.  We let $\vol(S^3 \setminus K')$ denote
\emph{simplicial volume}, i.e. the sum of the volumes of the
hyperbolic pieces in the geometric decomposition of $S^3 \setminus
K'$. 

We combine Theorem \ref{thm:periodic-smith} with a result of Agol,
Storm, and Thurston (see Theorem \ref{thm:alt-volume}) to give a bound
in terms of the diagram of $K'$.  We first make the following
definitions.

\begin{define}\label{def:twist}
For a knot or link $K$, we consider a diagram $D(K)$ as a 4--valent
graph in the plane, with over--under crossing information at each
vertex.  A link diagram $D$ is called \emph{prime} if any simple
closed curve that meets two edges of the diagram transversely bounds a
region of the diagram with no crossings.

Two crossings of a link diagram $D$ are defined to be equivalent if
there is a simple closed curve in the plane meeting $D$ in just those
crossings.  An equivalence class of crossings is defined to be a
\emph{twist region}.  The number of distinct equivalence classes is
defined to be the \emph{twist number} of the diagram, and is denoted
$\twist$.
\end{define}

Our definition of twist number agrees with that in
\cite{agol-storm-thurston}, and differs slightly from that in
\cite{fkp-07}.  The two definitions agree provided the diagram is
sufficiently reduced (i.e. \emph{twist reduced} in \cite{fkp-07}).  We
prefer Definition \ref{def:twist} as it does not require us to further
reduce diagrams.

\begin{corollary}\label{periodicalter}
With the notation and setting of Theorem \ref{thm:periodic-smith}
suppose, moreover, that $K'$ is alternating and hyperbolic, with prime
alternating diagram $D'$.  Then
$$ {\vol(S^3\setminus K)} \: \geq \: \left(1-\frac{2\sqrt 2
\pi^2}{p^2}\right)^{3/2} \, p\, v_8\, \left(\frac{\tw(D')}{2} -1
\right),
$$
where $v_8 =
	3.6638\dots$ is the volume of a regular ideal octahedron in
	$\HH^3$.\
\end{corollary}

By combining Theorem \ref{thm:periodic-smith} with recent results by
Agol \cite{agol:two-cusped} and Gabai, Meyerhoff, and Milley
\cite{gmm:smallest-cusped}, we obtain a universal estimate for the
volumes of periodic links. For ease of notation, define the function
$\psi\co \{ x\in \RR\, \co \, x{\geq 5.5}\} \to \RR$ by
$$	\psi(x) := \min \left\{2.828, \: 3.647 \left(1-\frac{2\sqrt 2 \,
	\pi^2}{x^2}\right)^{3/2} \right\}.
$$
Note that the right--hand term in the definition of $\psi$ is
greater than $2.828$ for $x \geq 14$.

\begin{theorem}\label{thm:universal-periodic}
Let $K$ be a hyperbolic periodic link in $S^3$, of period $p\geq 6$,
where we allow freely periodic links as well as those in which
the symmetry acts by rotation.  Then either
\begin{enumerate}
	\item $\vol(S^3 \setminus K) \: \geq \: p \cdot \psi(p)$, or
	\item $K$ is one of two explicit exceptions: a $5$--component link
of period $10$ whose quotient is $L(10,3) \setminus {\tt m003}$ or a $5$--component link
of period $15$ whose quotient is $L(15,4) \setminus {\tt m006}$. Here, ${\tt m003}$ and ${\tt m006}$ are manifolds in the SnapPea census; each of these two manifolds is the complement of a unique knot in the respective lens space.
\end{enumerate}
The estimate (1) is sharp for four freely periodic links, whose
periods are $14$, $18$, $19$, and $21$.
\end{theorem}

\subsection{Tangles and volumes} \label{sec:volume-results}
A \emph{tangle diagram} $T$ (or simply a \emph{tangle}) is a
graph contained in a unit square in the plane, with four 1--valent
vertices at the corners of the square, and all other vertices
4--valent in the interior.  Just as with knot diagrams, every
4--valent vertex of a tangle diagram comes equipped with over--under
crossing information.  Label the four 1--valent vertices as NW, NE,
SE, SW, positioned accordingly.

A tangle diagram is defined to be \emph{prime} if, for any simple
closed curve contained within the unit square which meets the diagram
transversely in two edges, the bounded interior of that curve contains
no crossings.
Two crossings in a tangle are equivalent if there is a simple closed
curve in the unit square meeting $D$ in just those crossings.
Equivalence classes are called \emph{twist regions}, and the number of
distinct classes is the \emph{twist number} of the tangle.

An alternating tangle is called \emph{positive} if the NE strand leads
to an over-crossing, and \emph{negative} if the NE strand leads to an
under-crossing.

The \emph{closure of a tangle} is defined to be the link diagram
obtained by connecting NW to NE and SW to SE by crossing--free arcs on
the exterior of the disk.
A \emph{tangle sum}, also called a \emph{Conway sum}, of tangles $T_1,
\dots, T_n$ is the closure of the tangle obtained by connecting
diagrams of the tangles $T_1, \dots, T_n$ linearly west to
east. Notice that if $T_1, \dots, T_n$ are all positive or all
negative, their tangle sum will be an alternating diagram.

Finally, we will call a tangle diagram $T$ an \emph{east--west twist}
if $\tw(T)=1$ and the diagram consists of a string of crossings
running from east to west.  The closure of such a diagram gives a
standard diagram of a $(2, q)$ torus link.

\begin{theorem} \label{thm:tangle-sum-new1}
Let $T_1, \dots, T_n$, $n\geq 12$, be tangles admitting prime,
alternating diagrams, none of which is an east--west twist.  Let $K$
be a knot or link which can be written as the Conway sum of the
tangles $T_1, \dots, T_n$, with diagram $D$.  Then $K$ is hyperbolic,
and
$$ \frac{v_8}{2}  \left(1 - \left( \frac{8\pi}{11.524 + n
	\sqrt[4]{2}} \right)^2 \right)^{3/2} \left(\tw(D) -3\right) \:
	\leq \: \vol(S^3-K) \: < \:10\, v_3 \, (\tw(D) - 1).$$ Here, $v_3 =
	1.0149\dots$ is the volume of a regular ideal tetrahedron and $v_8 =
	3.6638\dots$ is the volume of a regular ideal octahedron in
	$\HH^3$.\
\end{theorem}

The upper bound is due to Agol and D. Thurston
\cite{lackenby:alt-volume}.  The lower bound approaches
$(v_8/2)(\twist-3)$ as the number of tangles $n$ approaches infinity --
 similar to the (sharp) lower bound for alternating diagrams
proved by Agol, Storm, and Thurston \cite{agol-storm-thurston}.
However, Theorem \ref{thm:tangle-sum-new1} applies to more classes of
knots than alternating.  For example, it applies to large classes of
arborescent links (e.g. Montesinos links of length at least 12).  In
fact, our method of proof applies to links that are obtained by
summing up any number of ``admissible'' tangles, where the term
admissible includes, but is not limited to, alternating tangles,
tangles that admit diagrams containing at least seven crossings per
twist region and tangles whose closures are links of braid index $3$.

\subsection{Jones polynomial relations}
The volume conjecture of Kashaev and Murakami-Murakami asserts that
the volume of hyperbolic knots is determined by certain asymptotics of
the Jones polynomial and its relatives.  At the same time, recent
results \cite{dasbach-lin:volumeish, fkp-07} 
combined with a wealth of experimental evidence suggest a
coarse version of the volume conjecture: that the coefficients of the
Jones polynomial of a hyperbolic link should determine the volume of
its complement, up to bounded constants. To state the contribution of
the current paper to this coarse volume conjecture we need some
notation.  For a link $K$, we write its Jones polynomial in the form
$$J_K(t)= \alpha t^k+ \beta t^{k-1}+ \ldots + \beta' t^{m+1}+
\alpha' t^m, $$
so that the second and next-to-last coefficients of $J_K(t)$ are
$\beta$ and $\beta'$, respectively.  Dasbach and Lin proved
\cite{dasbach-lin:volumeish} that if $D(K)$ is a prime, alternating
diagram, then $\tw(D)=\abs{\beta} + \abs{\beta'}$. In
\cite{fkp-07}, we extended their results to give relations between the
coefficients of the Jones polynomial of links and the twist number of
link projections that contain at least three crossings per twist
region.  We further extend the result here.

Above, we defined the \emph{closure} of a tangle (also called the
\emph{numerator closure}) to be the link diagram obtained by
connecting NW to NE and SW to SE by simple arcs with no crossings.
The \emph{denominator closure} of the tangle is defined to be the
diagram obtained by connecting NW to SW, and NE to SE by simple arcs
with no crossings.  We say that a tangle diagram $T$ is \emph{strongly
alternating} if it is alternating and both the numerator and
denominator closures define prime diagrams.

\begin{theorem} \label{tangleJP}	
Let $T_1, \dots, T_n$ be alternating tangles whose Conway sum is a
knot $K$ with diagram $D$.  Define $T_+$ to be the result of joining
all the positive $T_i$ west to east, $T_-$ to be the result of joining
all the negative $T_i$ west to east.  Then, if both $T_+$ and $T_-$
are strongly alternating,
$$\frac{\twist}{2}-2 \; \leq \; \abs{\beta} + \abs{\beta'} \; \leq \;
  2 \, \twist.$$
\end{theorem}

If some $T_i$ is an east--west twist, then the denominator closure of
$T_+$ or $T_-$ will contain nugatory crossings, failing to be
prime.  Thus the hypotheses of Theorem \ref{tangleJP} imply that no
$T_i$ is an east--west twist. As a result, combining Theorems
\ref{thm:tangle-sum-new1} and \ref{tangleJP} gives

\begin{corollary}\label{tvolumish}
Let $K$ be a knot which can be written as the Conway sum of tangles
$T_1, \dots, T_n$. Let $T_+$ and $T_-$ be the sums of
the positive and negative $T_i$, respectively. Suppose that $n\geq
12$, and both $T_+$ and $T_-$ are strongly alternating. Then $K$ is
hyperbolic, and
$$ \frac{v_8}{4}\, \left(1 - \left( \frac{8\pi}{11.524 + n
\sqrt[4]{2}} \right)^2 \right)^{3/2} \left(\abs{\beta}+\abs{\beta'} -6
\right) \: \leq \: \vol(S^3-K) \:<\: 20 v_3 \, \left(\abs{\beta} +
\abs{\beta'} + \frac{3}{2}\right).
$$
\end{corollary}

The hypothesis that $K$ be a knot  is crucial in the statements of Theorem
\ref{tangleJP} and Corollary \ref{tvolumish}. Both
statements fail, for example, for the family of $(2, \cdots,2, -2, \cdots, -2)$
pretzel links. Theorem \ref{thm:tangle-sum-new1} implies the volume of $K$ will grow in an approximately linear fashion with the number of positive and negative $2$'s. On the other hand, using Lemma \ref{lem:oxs} below one can easily compute that $\abs{\beta} + \abs{\beta'} = 2$ for this family of links.

\subsection{Organization}
The proofs of our theorems bring together several very recent results
of Agol \cite{agol:two-cusped}; Agol, Storm, and Thurston
\cite{agol-storm-thurston}; Gabai, Meyerhoff, and Milley
\cite{gmm:smallest-cusped}; and the authors \cite{fkp-07}.  We survey
the results in Section \ref{sec:methods}.  In Section
\ref{sec:periodic}, we move on to periodic links to prove Theorems
\ref{thm:periodic-smith} and \ref{thm:universal-periodic} and
establish some corollaries.  Then, in Section \ref{sec:belt} we use
Adams' ``belted sum'' operations to study the behavior of hyperbolic
volume under the Conway summation of tangles, proving Theorem
\ref{thm:tangle-sum-new1}.  In Section \ref{sec:jones} we prove
Theorem \ref{tangleJP}.

\section{Recent estimates of hyperbolic volume and cusp area}\label{sec:methods}
In this section, we survey several recent results by Agol
\cite{agol:two-cusped}, Agol--Storm--Thurston
\cite{agol-storm-thurston}, the authors \cite{fkp-07}, and
Gabai--Meyerhoff--Milley \cite{gmm:smallest-cusped}, which we will
apply in later sections. Taken together, these theorems give powerful
structural results about the volumes of hyperbolic manifolds.
We also prove Theorem \ref{thm:gmm-cusp-estimate}, which follows
quickly from the above recent results, and will be important in
Section \ref{sec:belt}.

\subsection{Estimates from guts}

Let $M$ be a hyperbolic $3$--manifold, and $S \subset M$ an essential
surface. When we cut $M \setminus S$ along essential annuli, it
decomposes into a characteristic submanifold $B$ (the union of all
$I$--bundles in $M \setminus S$), and a hyperbolic component called
$\guts(M,S)$. Using Perelman's estimates for volume change under Ricci
flow with surgery, Agol, Storm, and Thurston proved the following
result.

\begin{theorem}[Theorem 9.1 of \cite{agol-storm-thurston}]\label{thm:guts}
Let $M$ be a finite--volume hyperbolic $3$--manifold, and let $S
\subset M$ be an essential surface. Then
$$\vol(M) \geq -v_8 \, \chi(\guts(M,S)).$$
\end{theorem}

Combining Theorem \ref{thm:guts} with Lackenby's analysis of
checkerboard surfaces in alternating link complements
\cite{lackenby:alt-volume} gives the following result, which bounds volume
based on diagrammatic properties.

\begin{theorem}[Corollary 2.2 of \cite{agol-storm-thurston}]\label{thm:alt-volume}
Let $D(K)$ be a prime, alternating link diagram with
$\twist \geq 2$. Then $K$ is hyperbolic, and
$$\vol(S^3 \setminus K) \geq \frac{v_8}{2} \, (\twist - 2).$$
\end{theorem}

More recently, Agol showed that every two--cusped hyperbolic
$3$--manifold contains an essential surface with non-trivial guts
\cite{agol:two-cusped}. Using Theorem \ref{thm:guts}, he obtained

\begin{theorem}[Theorem 3.4 of \cite{agol:two-cusped}]\label{thm:two-cusped-volume}
Let $M$ be an orientable hyperbolic $3$--manifold with two or more cusps. Then
$$\vol(M) \geq v_8,$$
with equality if and only if $M$ is the
complement of the Whitehead link or its sister ({\tt m129} or {\tt
m125} in the notation of the SnapPea census).
\end{theorem}

\subsection{Bounding volume change under Dehn filling}
Given a 3--manifold $M$ with at least $k$ torus boundary components,
we use the following standard terminology.  For the $i$-th torus
$T_i$, let $s_i$ be a \emph{slope} on $T_i$, that is, an isotopy class
of simple closed curves.  Let $M(s_1, \dots, s_k)$ denote the manifold
obtained by Dehn filling along the slopes $s_1$, \dots, $s_k$.

When $M$ is hyperbolic, each torus boundary component of $M$
corresponds to a cusp.  Taking a maximal disjoint horoball
neighborhood about each of the cusps, each torus $T_i$ inherits a Euclidean
structure, well--defined up to similarity.  The slope $s_i$ can then
be given a geodesic representative.  We define the \emph{slope length}
of $s_i$ to be the length of this geodesic representative.  Note that
when $k>1$, this definition of slope length depends on the choice of
maximal horoball neighborhood.  The authors recently showed the
following result.

\begin{theorem}[Theorem 1.1 of \cite{fkp-07}]\label{thm:fkp-dehn-filling}
  Let $M$ be a complete, finite--volume hyperbolic manifold with
  cusps.  Suppose $C_1, \dots, C_k$ are disjoint horoball
  neighborhoods of some subset of the cusps.  Let $s_1, \dots, s_k$ be
  slopes on $\partial C_1, \dots, \partial C_k$, each with length
  greater than $2\pi$.  Denote the minimal slope length by $\lmin$.
  Then $M(s_1, \dots, s_k)$ is a hyperbolic manifold,
  and
  $$ \vol(M(s_1, \dots, s_k)) \: \geq \:
  \left(1-\left(\frac{2\pi}{\lmin}\right)^2\right)^{3/2} \vol(M).$$
\end{theorem}

\subsection{Mom technology}
In a series of recent papers \cite{gmm:mom-tech, gmm:smallest-cusped, milley:smallest-volume},
Gabai, Meyerhoff, and Milley developed the theory of \emph{Mom
manifolds}. A \emph{Mom-n structure} $(M,T,\Delta)$ consists of a
compact 3--manifold $M$ whose boundary is a union of tori, a preferred
boundary component $T$, and a handle decomposition $\Delta$ of the
following type.  Starting from $T\times I$, $n$ 1--handles and $n$
2--handles are attached to $T\times 1$ such that each 2--handle goes
over exactly three 1--handles, counted with multiplicity.  Furthermore,
each 1--handle encounters at least two 2--handles, counted with
multiplicity. We say that $M$ is a \emph{Mom-n} if it possesses a
Mom-$n$ structure $(M,T,\Delta)$.

In \cite{gmm:mom-tech}, Gabai, Meyerhoff, and Milley enumerated all
the hyperbolic Mom-2's and Mom-3's (there are 21 such manifolds
in total). In \cite{gmm:smallest-cusped}, they showed that every
cusped hyperbolic manifold of sufficiently small volume (or cusp area)
must be obtained by Dehn filling a Mom-2 or Mom-3 manifold:

\begin{theorem}[\cite{gmm:smallest-cusped}]\label{thm:internal-mom}
Let $M$ be a cusped, orientable hyperbolic $3$--manifold. Assume that
$\vol(M) \leq 2.848$ or that a maximal horoball neighborhood $C$ of
one of its cusps has $\area(\bdy C) \leq 3.78$. Then $M$ is obtained
by Dehn filling on one of the 21 Mom-2 or Mom-3 manifolds.
\end{theorem}

\begin{proof}
The volume part of the theorem is explicitly stated as Theorem 1.1 of
\cite{gmm:smallest-cusped}. The cusp area part of the statement
follows by evaluating Gabai, Meyerhoff, and Milley's cusp area
estimates \cite[Lemmas 4.6, 4.8, 5.4, 5.6, and
5.7]{gmm:smallest-cusped} on the parameter space of all ortholengths 
corresponding to manifolds
without a Mom-2 or Mom-3 structure. The rigorous C++ and Mathematica code to
construct and evaluate those estimates was helpfully supplied by
Milley \cite{milley:computation}.
\end{proof}

Because each of the Mom-2 and Mom-3 manifolds has volume significantly
higher than 2.848, Theorem \ref{thm:fkp-dehn-filling} bounds the
length of the slope along which one must fill a Mom manifold to obtain
$M$.  Thus,
Theorem \ref{thm:internal-mom} combined with Theorem \ref{thm:fkp-dehn-filling}
reduces
the search for small--volume manifolds
to finitely many Dehn fillings of the 21 Mom-2's and Mom-3's.

\begin{corollary}[Theorem 1.2 of \cite{milley:smallest-volume}]\label{cor:small-vol-cusped}
Let $M$ be a cusped, orientable hyperbolic manifold whose volume is at
most $2.848$. Then $M$ is one of the SnapPea census manifolds {\tt
m003}, {\tt m004}, {\tt m006}, {\tt m007}, {\tt m009}, {\tt m010},
{\tt m011}, {\tt m015}, {\tt m016}, or {\tt m017}. In particular,
every cusped hyperbolic manifold with $\vol(M) \leq 2.848$ can be
obtained by Dehn filling 
two cusps of the $3$--chain link complement in Figure \ref{fig:chain-link}.
\end{corollary}

Theorem \ref{thm:internal-mom} can also be employed to give universal
estimates for the cusp area of those manifolds that have two or more
cusps:

\begin{figure}
\begin{center}
\includegraphics{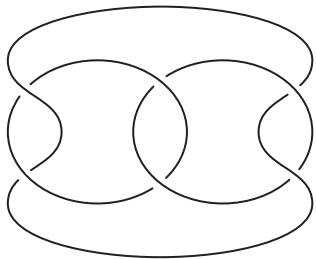}
\caption{The complement of the $3$--chain link is the only Mom-2 or
Mom-3 manifold with more than two cusps.}
\label{fig:chain-link}
\end{center}
\end{figure}

\begin{theorem}\label{thm:gmm-cusp-estimate}
Let $M$ be an orientable hyperbolic $3$--manifold with two or more
cusps. Suppose that $M$ contains a belt (an essential twice--punctured
disk). If $C$ is a maximal neighborhood of one of the cusps of $M$, then
$$\mathrm{area}(\bdy C) \geq 3.78.$$
\end{theorem}

\begin{remark}
The hypothesis that $M$ contains a belt should be
unnecessary. However, proving the theorem without this hypothesis
would require studying infinitely many fillings of the $3$--chain link
in Figure \ref{fig:chain-link}.
\end{remark}

\begin{proof}
Theorem \ref{thm:internal-mom} implies that every cusped hyperbolic
manifold either has cusp area at least $3.78$, or is obtained by Dehn
filling on one of the Mom-2 or Mom-3 manifolds. Among these 21 Mom
manifolds, 20 have exactly two cusps. Thus, if $M$ is obtained by filling
on one of these 20 manifolds, the filling must be trivial and it already \emph{is} one of the Mom
manifolds. Individual verification shows that a maximal neighborhood
of any cusp of any of the Mom-2 or Mom-3 manifolds has area at least
$4$ (with the minimum of $4$ realized by the Whitehead
link). Therefore, $M$ either has cusp area at least $3.78$, or is
obtained by Dehn filling one cusp of the single $3$--cusped Mom
manifold $N$, namely the complement of the $3$--chain link depicted in
Figure \ref{fig:chain-link}.

\begin{prop}\label{prop:3-chain}
Let $M$ be a hyperbolic $3$--manifold obtained by filling one cusp of
the \linebreak $3$--chain link complement $N$. Suppose that $M$
contains an essential twice--punctured disk. If $C$ is a maximal
neighborhood of one of the cusps of $M$, then $\mathrm{area}(\bdy C)
\geq 4$.
\end{prop}

\begin{proof}
Suppose that $M$ contains an essential twice--punctured disk
$P$. Isotope $P$ to minimize its intersection number with the core of
the solid torus added during Dehn filling. Then $S = P \cap N$ is an
essential surface in $N$; more precisely, it is an essential sphere
with $(n+3)$ holes, where $n$ of its boundary circles run in parallel
along the filling slope. Since every thrice--punctured sphere in $N$
meets all three cusps (and thus becomes an essential annulus after
filling along one of its boundary circles), we can conclude that $n
\geq 1$.

Now, expand a maximal horospherical neighborhood $H$ of the cusp of
$N$ that we are filling. Consider the length $\ell$ of the filling
slope along $\bdy H$. Since $S \cap \bdy H$ consists of $n$ distinct
circles of that slope, a result of Agol and Lackenby (see
\cite[Theorem 5.1]{agol:6theorem} or \cite[Lemma
3.3]{lackenby:surgery}) implies that the total length of those circles
is
$$n \, \ell \: \leq \: -6 \, \chi(S) \: = \: 6(n+1) \: \leq \: 12 n.$$
Therefore, $M$ is obtained by filling one cusp of $N$ along a slope of
length at most $12$.

To complete the proof, we enumerate the slopes that have length at
most $12$. Note that since the symmetry group of $N$ permutes all
three cusps, it suffices to consider a single cusp. In complex
coordinates on this maximal cusp, the knot--theoretic longitude is a
translation by $4$, while the meridian is a translation by
$\frac{3}{2} + \frac{\sqrt{7}}{2} i$. Thus the slopes on a cusp of $N$
that have length at most $12$ are:
\begin{equation}\label{eq:chain-short-slopes}
\begin{array}{c}
1/0 \\ 
-7 \quad  -\!6 \quad \cdots \quad  3 \quad \: 4 \\
-7/2 \quad -\!5/2 \quad  -\!3/2 \quad  -\!1/2 \quad \: 1/2 \\
-8/3 \quad -\!7/3 \quad  -\!5/3  \quad  -\!4/3 \quad  -\!2/3 \quad -\!1/3 \\ 
-7/4 \quad -\!5/4
\end{array}
\end{equation}

Martelli and Petronio \cite{martelli-petronio:magic-manifold} have
shown that the non-hyperbolic fillings of one cusp of $N$ are exactly
the fillings along slope $\infty, -3, -2, -1, \, 0$. 
For each of the $21$ remaining slopes, SnapPea finds (an approximate solution for) a
hyperbolic structure on the filled manifold. H.\ Moser's thesis \cite{moser:thesis} then implies that the true hyperbolic structure on each of these manifolds is indeed $\varepsilon$--close to the one found by SnapPea. In each case, the cusp
area is bounded below by $4$.\end{proof}

Proposition \ref{prop:3-chain} completes the proof of Theorem
\ref{thm:gmm-cusp-estimate}.\end{proof}

As a closing remark, we point out that
among the hyperbolic fillings of the $3$--chain link listed in
(\ref{eq:chain-short-slopes}), only the Whitehead link complement
contains a belt. In other words, a topological analysis of these
manifolds shows that the Whitehead link is the only manifold
satisfying the hypotheses of Proposition \ref{prop:3-chain}. 
Since we do not need this stronger statement in the sequel, we omit
the details.

\section{Volume estimates for periodic links}\label{sec:periodic}
Let $K$ be a periodic link and let $h: {S}^3 \rightarrow {S}^3$ be an
orientation preserving diffeomorphism of order $p$ with $h(K)=K$, such that the set
of fixed points $C_h$ of $h$ is a circle that is either
disjoint from $K$ or is a component of $K$.  By Smith theory and the solution to
the Smith conjecture \cite{smith}, $C_h$ is the trivial knot and $h$
is conjugate to a rotation with axis $C_h$. The quotient of the action
of $h$ on $K$ is a link $K'$, called the \emph{quotient} of $K$.  Let $C'_h$
denote the quotient of the axis $C_h$ under the action of $h$ on
$S^3$.

\begin{theorem}\label{periodic1}
Let $K$ be a periodic hyperbolic link in $S^3$ of period $p\geq
6$. Let $C'_h$ be the quotient of the fixed point set under $h$ and
let $K'$ be the quotient link of $K$.  Then, $L_h:=K'\cup C'_h$ is a
hyperbolic link, and
$$ p \; \left(1-\frac{2\sqrt 2\pi^2}{p^2}\right)^{3/2}
  \vol(S^3\setminus L_h) \: \leq \: \vol(S^3\setminus K) \: \leq \:
  p\;\vol(S^3\setminus L_h) .
$$
\end{theorem}

\begin{proof}
The Mostow--Prasad rigidity theorem implies that $h$ can be homotoped
to a hyperbolic isometry $h\co S^3\setminus K \rightarrow
S^3\setminus K$.  Since $S^3\setminus K$ is a Haken 3--manifold, a
result of Waldhausen \cite{waldhausen:suff-large} implies that $h$ can
actually be isotoped to a hyperbolic isometry. Thus $C_h$ is either a component of $K$, or else it is
a closed geodesic in $S^3 \setminus K$.  It follows that $S^3 \setminus (K \cup
C_h)$ is hyperbolic (in the case that $C_h$ is a component of $K$ we take $K \cup
C_h=K$).  Now the quotient of the action $h\co S^3\setminus
(K\cup C_h) \rightarrow S^3\setminus (K\cup C_h)$, which is
$S^3\setminus L_h$, is also hyperbolic.  The quotient map
$$S^3\setminus (K\cup C_h) \longrightarrow S^3\setminus L_h$$
is a covering of degree $p$.  Thus
$$\vol(S^3\setminus (K\cup C_h) )=p \, \vol(S^3\setminus L_h ).$$

If $C_h$ is a component of $K$, we are done. Otherwise, $S^3\setminus K$ is obtained from $S^3\setminus (K\cup C_h)$ by Dehn filling $C_h$ along the
meridian $m$.  This meridian covers the meridian $m'$ of
$C'_h$ $p$ times.  By work of Adams \cite{adams:waist2}, the length
of $m'$ satisfies $l(m')\geq \sqrt[4]{2}$.  Thus, $l(m)\geq p\,
\sqrt[4]{2}$.  For $p\geq 6$ we have $l(m)\geq p\,
\sqrt[4]{2}>2\pi$.  Now Theorem \ref{thm:fkp-dehn-filling} applies,
and we conclude
\begin{eqnarray*}
\vol(S^3\setminus K) &\geq&
\left(1-\frac{2\sqrt 2\pi^2}{p^2}\right)^{3/2}\vol(S^3\setminus (C_h\cup K))\\
 &=&
 \left(1-\frac{2\sqrt 2\pi^2}{p^2}\right)^{3/2}\, p\, \vol(S^3\setminus L_h).\\
\end{eqnarray*}

As for the upper bound, we note that volume strictly decreases under
Dehn filling \cite[Corollary 6.5.2]{thurston:notes}. Thus, if $C_h$ is
not already a component of $K$, we have
$$p\;\vol(S^3\setminus L_h) =\vol(S^3\setminus (K\cup C_h)
)>\vol(S^3\setminus K).$$

\vspace{-3ex}
\end{proof}

Next we derive Theorem \ref{thm:periodic-smith} from Theorem
\ref{periodic1}: To that end, for a $3$--manifold $M$ we will let
$\left\| M \right\|$ denote the \emph{Gromov norm} of $M$. By
\cite[Theorem 6.5.4]{thurston:notes}, if $M$ is hyperbolic then
$\vol(M)=v_3 \left\| M \right\|$. More generally, $v_3 \left\| M
\right\|$ is the simplicial volume of $M$, equal to the sum of volumes
of the hyperbolic pieces in the geometric decomposition of $M$.

\begin{proof}[Proof of Theorem \ref{thm:periodic-smith}] 
If the axis $C_h$ is not already a component of $K$, the complement
$S^3\setminus K'$ is obtained by Dehn filling from $S^3\setminus
L_h$. We note that $K'$ need not be hyperbolic.  By \cite[Corollary
6.5.2]{thurston:notes}, we have $\left\|S^3\setminus
L_h\right\|>\left\|S^3\setminus K'\right\|$.  Since, by Theorem
\ref{periodic1}, $L_h$ is hyperbolic, $\vol(S^3\setminus L_h)=v_3
\left\|S^3\setminus L_h\right\|$.  Combining these facts with the
left-hand inequality of Theorem \ref{periodic1} gives
$$\vol(S^3\setminus K) \geq \left(1-\frac{2\sqrt
	2\pi^2}{p^2}\right)^{3/2}\, p\,v_3 \,  \left\|S^3\setminus
K'\right\|.$$ 

\vspace{-4.5ex}
\end{proof}

Now, we turn our attention to Theorem
\ref{thm:universal-periodic}.  Define $\psi\co \{ x\in \RR\, \co \, x{\geq 5.5}\} \to \RR$ by
$$	\psi(x) := \min \left\{2.828, \: 3.647 \left(1-\frac{2\sqrt 2 \,
	\pi^2}{x^2}\right)^{3/2} \right\}.
$$

\begin{theorem}[Theorem \ref{thm:universal-periodic}]
Let $K$ be a hyperbolic periodic link in $S^3$, of period $p\geq
6$.
Then either
\begin{enumerate}
\item $\vol(S^3 \setminus K) \: \geq \: p \cdot \psi(p)$, or
\item $K$ is one of two explicit exceptions: a $5$--component link
of period $10$ whose quotient is $L(10,3) \setminus {\tt m003}$ or a $5$--component link
of period $15$ whose quotient is $L(15,4) \setminus {\tt m006}$. 
\end{enumerate}
Estimate (1) is sharp for four freely periodic links,
whose periods are $14$, $18$, $19$, and $21$.
\end{theorem}

\begin{proof}
Let $h:S^3 \to S^3$ be the diffeomorphism of order $p$ that sends $K$
to itself. As discussed in the introduction, the solutions to the
Smith conjecture \cite{smith} and the spherical spaceforms conjecture
\cite{morgan-tian} imply that we may take $h$ to be an orthogonal
action by an element of $SO(4)$. We need to consider two cases:
either $h$ fixes an invariant axis $C_h$, or $h^i$ acts on $S^3$ without fixed
points, for all $0<i<p$.

If $h$ has an invariant axis $C_h$, then Theorem \ref{periodic1}
applies, and
$$\vol(S^3\setminus K) \: \geq \: p \; \left(1-\frac{2\sqrt
2\pi^2}{p^2}\right)^{3/2} \vol(S^3\setminus L_h).$$
Now, because $L_h$
is a hyperbolic link of two or more components, Agol's Theorem
\ref{thm:two-cusped-volume} gives $\vol(S^3\setminus L_h) \geq 3.663$,
completing the argument in this case.

If  $h^i$ acts on $S^3$ without fixed points, for all $0<i<p$, the quotient of $S^3$ is a
lens space $L(p,q)$ and the quotient of $S^3 \setminus K$ is a
hyperbolic manifold $M$, obtained as the complement of a link in
$L(p,q)$. Thus
$$\vol(S^3 \setminus K) = p \cdot \vol(M).$$
If $\vol(M) \geq 2.828$, then $K$ satisfies the statement of the
theorem. On the other hand, if $\vol(M) \leq 2.848$, then $M$ is one
of the ten one--cusped manifolds listed in Corollary
\ref{cor:small-vol-cusped}. Thus, to complete the proof, it suffices
to enumerate all of the ways in which each of these ten manifolds
occurs as the complement of a knot in a lens space. Because each
manifold in Corollary \ref{cor:small-vol-cusped} is a filling of two
cusps of the complement $N$ of the $3$--chain link of Figure 1, we can use the extensive tables
compiled by Martelli and Petronio \cite[Section
A.1]{martelli-petronio:magic-manifold} to enumerate their lens space
fillings:

\begin{center}
\begin{tabular}{|c|c|c|c|c|}
\hline
Manifold & Alternate name & Volume & Surgery on $N$ & Lens space fillings \\
\hline
{\tt m003} & figure--8 sister & 2.0298... & $N(1, -4)$ & $L(5,1), \, L(10, 3)$ \\
{\tt m004} & figure--8 knot & 2.0298... & $N(1, 2)$ & $S^3$ \\
{\tt m006} & & 2.5689... & $N(1, -3/2)$ & $L(5,2), \, L(15,4)$ \\
{\tt m007} & & 2.5689... & $N(1, -1/2)$ &  $L(3,1)$ \\
{\tt m009} & p. torus bundle $LLR$ & 2.6667... & $N(1, 3)$ & $L(2,1)$ \\
{\tt m010} & p. torus bundle $-LLR$ & 2.6667... & $N(1, -5)$ & $L(6,1)$ \\
{\tt m011} & & 2.7818... & $N(-3/2, -5)$ & $L(9, 2), \, L(13,4)$ \\
{\tt m015} & $5_2$ knot & 2.8281... & $N(1, 1/2)$ & $S^3$ \\
{\tt m016} & $(-2,3,7)$ pretzel knot & 2.8281... & $N(-3/2, -1/2)$ & $S^3, \, L(18,5), \, L(19,7)$ \\
{\tt m017} &  & 2.8281... & $N(1, -5/2)$ & $L(7,2), \, L(14,3), \, L(21,8)$ \\
\hline
\end{tabular}
\end{center}

The proof will be complete after several observations. First, we may
ignore lens spaces $L(p,q)$ with $p \leq 5$, because we have assumed
$p \geq 6$. Second, the two exceptions to the theorem are obtained by
lifting to $S^3$ the knots $(L(10,3) \setminus {\tt m003})$ and
$(L(15,4) \setminus {\tt m006})$. Homology considerations
show that both of these exceptions are $5$--component links. Third,
even though $L(6,1)$, $L(9,2)$, and $L(13,4)$ are obtained by filling manifolds of volume less
than $2.828$, the corresponding links satisfy the theorem
because $3.647 \left(1- 2\sqrt 2 \, \pi^2 /6^2\right)^{3/2} < 2.666$ and
$$
3.647 \left(1- \frac{2\sqrt 2 \, \pi^2}{9^2} \right)^{3/2} \:<\:\: 
3.647 \left(1- \frac{2\sqrt 2 \, \pi^2}{13^2} \right)^{3/2} \:<\:\:
2.7818.$$
Finally, the four examples demonstrating the sharpness of the
theorem are the 18--fold and 19--fold covers of ${\tt m016}$ and the
14--fold and 21--fold covers of ${\tt m017}$.
\end{proof}

Note that if the link $K$ in Theorem \ref{thm:universal-periodic} is not
freely periodic, then the volume is actually bounded by the quantity
on the right in the definition of $\psi(n)$.

\section{Belted sums and Conway sums}\label{sec:belt}

\subsection{Belted sums}

\begin{figure}[h]
\input{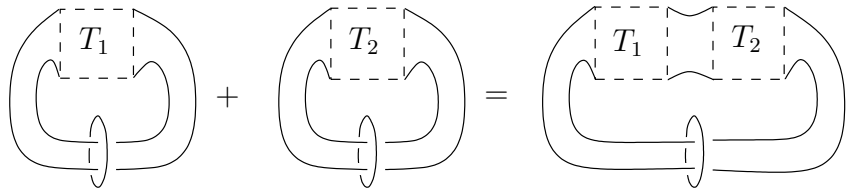}
\caption{Belted sum.}
\label{fig:beltsum}
\end{figure}

Let $T$ be a tangle diagram. Given $T$, we may form a link diagram as
follows.  First, form the closure of $T$ by connecting NE to NW, and
SE to SW. Then, add an extra component $C$ that lies in a plane
orthogonal to the projection plane and encircles the two unknotted
arcs that we have just added to $T$. See the left of Figure
\ref{fig:beltsum}.  We call the resulting link the \emph{belted tangle
corresponding to $T$}, or simply a \emph{belted tangle}.  Note that
$C$ bounds a 2--punctured disk $S$ in the complement of the link.  We
will call the link component $C$ the \emph{belt component} of the
link.  We will only be interested in belted tangles admitting
hyperbolic structures.

Given two hyperbolic belted tangles corresponding to $T_1$ and $T_2$,
with complements $M_1$ and $M_2$, belt components $C_1$ and $C_2$, and
2--punctured disks $S_1$ and $S_2$, we form the complement of a new
belted tangle as follows.  Cut each manifold $M_i$ along the surface
$S_i$, and then glue two manifolds with two 2--punctured
disks as boundary.  Since there is a unique hyperbolic structure on a
2--punctured disk we may glue $M_1$ to $M_2$ by an isometry that glues $C_1$ to
$C_2$.  The result is the complement of a new belted tangle.  See
Figure \ref{fig:beltsum}.  We call this new belted tangle the
\emph{belted sum} of the tangles $T_1$ and $T_2$.  Belted sums were
studied extensively by Adams \cite{adams:3-punct}.
Note that the Conway sum of $T_1$ and $T_2$ is obtained by meridional
Dehn filling on the belt component of the belted sum of $T_1$ and
$T_2$.

\subsection{Arc lengths on belted tangles}

Consider a maximal neighborhood $C$ of the cusp corresponding to the
belt component.  Denote by the \emph{width} the length of the shortest
nontrivial geodesic arc running from the 2--punctured disk to itself
on $\partial C$.  Adams \emph{et al} observed that the length of the
shortest nontrivial arc from an embedded totally geodesic surface to
itself is bounded below by $1$ (see \cite[Theorem 4.2]{adams:II} or
\cite[Theorem 1.5]{adams:quasi-fuch}).  In the case at hand, their
result gives the following.

\begin{lemma}
	The width of a belt component of a belted tangle is at least $1$.
\label{lemma:meridian-1}
\end{lemma}

Note that since the 2--punctured disk intersects the cusp in a
longitude, the meridian must be at least as long as the
width.  We will also need bounds on the length of a longitude.

\begin{lemma}
The length of the longitude of a belt component is at most $4$, and at
least $\sqrt[4]{2}$.
\label{lemma:long-4}
\end{lemma}

\begin{proof}
Both bounds are due to Adams. In 
\cite{adams:waist2}, he proves that if $M$ is not the
complement of the figure--8 knot or the $5_2$ knot, then the shortest
curve has length at least $\sqrt[4]{2}$.

As for the upper bound, the length of the longitude is maximal when
the maximal cusp in $M$ restricts to a maximal cusp on the 3--punctured sphere.
By \cite[Theorem 2.1]{adams:cusps-collars-systoles}, the length of a
maximal cusp on the 3--punctured sphere is $4$.
\end{proof}

We need to determine a maximal cusp corresponding to the belt
component of a belted sum of two tangles, $T_1$ and $T_2$.  When we
expand a horoball neighborhood about this cusp, the cusp neighborhood  may bump itself in
one component of the belt sum before it bumps in the other.  When the
cusp bumps itself, it determines a longitude of the belt component.
Thus the longitude of the belt component of the belted sum will have
length equal to the minimum of the longitude lengths of $T_1$ and
$T_2$.  Say this minimum occurs in $T_1$.  Then the length of any arc
running from 3--punctured sphere to 3--punctured sphere in $T_2$ will
be scaled by the ratio of the length of the longitude of $T_1$ and the
length of the longitude of $T_2$.
In particular, the width of the belted sum will not necessarily be the
width of $T_1$ plus the width of $T_2$, but rather the width of $T_1$
plus the width of $T_2$ times the ratio of the longitude length of
$T_1$ to the longitude length of $T_2$.


\begin{lemma}\label{lemma:general-width}
Let $T$ be a belted tangle obtained as the belted sum of $n$
hyperbolic belted tangles $T_1, \dots, T_n$.  Let $\ell$ be the length
of the shortest longitude of a belt component of the $T_j$.  Then the
width of the belt component of $T$ is at least
$$w \geq \frac{3.78}{\ell} + (n-1)\, \frac{\ell}{4}.$$
\end{lemma}

\begin{proof}
Without loss of generality, suppose $T_1$ has the shortest longitude.
By Theorem \ref{thm:gmm-cusp-estimate}, the cusp area corresponding to
the belt component is at least $3.78$.  Thus the width of $T_1$ is at
least $3.78/\ell$.  By Lemma \ref{lemma:long-4}, the longitudes of the
other $T_j$'s are at most $4$, and by Lemma \ref{lemma:meridian-1},
the widths of these are at least $1$.  When we do the belted sum, the
longitudes will rescale to be length $\ell$, and the widths will
rescale to be at least $\ell/4$.  Thus the total width will be at
least $w \geq 3.78/\ell + (n-1)(\ell/4).$
\end{proof}

\subsection{Volumes and belted tangles}

\begin{lemma}\label{lemma:vol-alt}
Let $T$ be a prime, alternating tangle that is not an east--west
twist.  Let $L$ denote the belted tangle corresponding to $T$.  
Then $L$ is hyperbolic.
Furthermore, 

\noindent
(A) If $1/n$ Dehn filling along the belt component adds a new
twist region to the closure of $T$, then
$$\vol(S^3 \setminus L) \geq \frac{v_8}{2}\,(\tw(T)-1).$$

\noindent
(B) If $1/n$ Dehn filling along the belt component adds crossings
to an existing twist region in the closure of $T$, then
$$\vol(S^3 \setminus L) \geq \frac{v_8}{2}\,(\tw(T)-2).$$
\end{lemma}

\begin{proof}
Let $L(n)$ denote the link formed by performing $1/n$ Dehn filling on
the belt component of $L$, where $n$ is positive or negative depending
on which sign makes $L(n)$ alternating.  When we form $L(n)$, we may
either add a new twist region to the closure of $T$, or we may add
additional crossings to an existing twist region.  In either
case the link $L(n)$ has at least two twist regions, since $T$ is not
an east--west twist, hence it is hyperbolic.

In case (A), Theorem \ref{thm:alt-volume} implies the volume of $S^3 \setminus L(n)$ is at
least $v_8/2( (\tw(T)+1) -2 )$.  In case (B), Theorem
\ref{thm:alt-volume} implies the volume of $S^3 \setminus L(n)$  is at least $v_8/2( \tw(T)
-2)$.  Because volume goes down under Dehn filling, these lower bounds
on the volume of $\vol(S^3 \setminus L(n))$ also serve as lower bounds
on $\vol(S^3 \setminus L)$.
\end{proof}

\begin{lemma}\label{lemma:belt-volume}
Let $T_1, \ldots, T_n$ prime, alternating tangle diagrams, none of
which is an east--west twist. Let $D(K)$ be the Conway sum of $T_1,
\ldots, T_n$, and let $L$ be the belted sum of these tangles. Then
$$\vol(S^3 \setminus L) \geq \frac{v_8}{2}\,(\tw(D)-3).$$
\end{lemma}

\begin{proof}
Because we formed the belted sum $L$ by gluing belted tangles along
totally geodesic $2$--punctured disks, the volume of $L$ will remain
unchanged if we permute the order of the $T_i$. Thus, without loss of
generality, we may assume that $T_1, \ldots, T_r$ are positive tangles
and $T_{r+1}, \ldots, T_n$ are negative tangles. Furthermore, if the
$T_i$ are all positive or all negative, then $D(K)$ is a prime,
alternating diagram, and the result follows by Lemma \ref{lemma:vol-alt}. 
Thus we may assume that $0<r<n$.

With these assumptions, let $D_+$ be the Conway sum and $L_+$ be the
belted sum of $T_1, \ldots, T_r$. Let $D_-$ be the Conway sum and
$L_-$ be the belted sum of $T_{r+1}, \ldots T_n$. Then each of $D_+$
and $D_-$ is the closure of a prime, alternating tangle. Thus, by
Lemma \ref{lemma:vol-alt},
$$\vol(S^3 \setminus L_+) \geq \frac{v_8}{2}\,(\tw(D_+)-2) \qquad
\mbox{and} \qquad 
\vol(S^3 \setminus L_-) \geq \frac{v_8}{2}\,(\tw(D_-)-2),
$$
with a sharper estimate if either $D_+$ or $D_-$ falls into case (A)
of the Lemma.

Suppose that either $D_+$ or $D_-$ falls into case (A) of Lemma
\ref{lemma:vol-alt}. Then, since equivalent crossings remain
equivalent after gluing, we have $\tw(D_+) + \tw(D_-) \geq \twist$,
and thus
$$\vol(S^3 \setminus L) \: = \: \vol(S^3 \setminus L_+) + \vol(S^3
\setminus L_-) \: \geq \: \frac{v_8}{2}\,(\tw(D)-3).$$
On the other hand, suppose that both $D_+$ and $D_-$ fall into case
(B) of Lemma \ref{lemma:vol-alt}. Then $1/n$ Dehn filling along the
belt component of both $L_+$ and $L_-$ adds crossings to existing
twist regions of both $D_+$ and $D_-$. In this situation, the
crossings in these two twist regions become equivalent when we join
$D_+$ and $D_-$. Thus $\tw(D_+) + \tw(D_-) \geq \twist + 1$, and
$$\vol(S^3 \setminus L) \: \geq \: \frac{v_8}{2}\, (\tw(D_+) +
\tw(D_-) - 4) \: \geq \: \frac{v_8}{2}\,(\tw(D)-3).$$

\vspace{-3ex}
\end{proof}

We may now prove Theorem \ref{thm:tangle-sum-new1}, which was stated
in the Introduction.

\begin{proof}[Proof of Theorem \ref{thm:tangle-sum-new1}]

Let $L$ be the belted sum of $T_1$, $\dots$, $T_n$.  We obtain $K$ by
meridional filling on the belt component of $L$.  By Lemma
\ref{lemma:belt-volume}, $\vol(S^3 \setminus L) \geq
v_8/2\,(\tw(D)-3)$.  Thus, using Theorem \ref{thm:fkp-dehn-filling},
we can estimate the volume of $S^3 \setminus K$ once we estimate the
meridian length of the belt.  To apply Theorem
\ref{thm:fkp-dehn-filling}, we also need to ensure that this length is
at least $2\pi$.  The meridian is at least as long as the width, which
by Lemma \ref{lemma:general-width} is at least $3.78/\ell +
(n-1)(\ell/4)$.

By Lemma \ref{lemma:long-4}, $\ell \in [\sqrt[4]{2}, 4]$.  Thus we
need to minimize the quantity
$$3.78/\ell + (n-1)(\ell/4)$$
over the interval $[\sqrt[4]{2}, 4]$.  For $n\geq 12$, we find this is
an increasing function of $\ell$, so the minimum value occurs when
$\ell = \sqrt[4]{2}$.  Hence the meridian will have length at least
$$
\lmin \: \geq \: \frac{3.78}{\sqrt[4]{2}} +
(n-1)\frac{\sqrt[4]{2}}{4} \: > \: \frac{11.524 + n \sqrt[4]{2}}{4} \, ,
$$
which is greater than $2\pi$ for $n \geq 12$. Thus Theorem
\ref{thm:fkp-dehn-filling} applies, and we obtain
\begin{eqnarray*}
\vol(S^3 \setminus K) &\geq&
\left(1-\left(\frac{2\pi}{\ell_{min}}\right)^2\right)^{3/2}\vol(S^3 \setminus L)\\
 &\geq& \left(1 -
	  \left( \frac{8\pi}{11.524 + n \sqrt[4]{2}} \right)^2
	\right)^{3/2} 
\frac{v_8}{2}\, \left( \tw(D) - 3\right).
\end{eqnarray*}

\vspace{-2ex}
\end{proof}

\section{The Jones polynomial and tangle addition}\label{sec:jones}
In this section,
we will prove Theorem \ref{tangleJP}, which gives Corollary
\ref{tvolumish}.

\subsection{ Adequate link preliminaries}
\begin{figure}[h]
	\centerline{\input{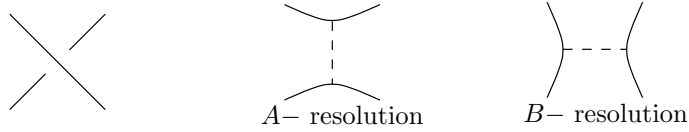}}
\caption{$A$-- and $B$--resolutions of a crossing.}
\label{fig:splicing}
\end{figure}

We begin by recalling 
some terminology and notation from \cite{dasbach-futer...} and \cite{fkp-07}.  
Let $D$ be a link diagram, and
$x$ a crossing of $D$.  Associated to $D$ and $x$ are two link
diagrams, each with one fewer crossing than $D$, called the
\emph{$A$--resolution} and \emph{$B$--resolution} of the crossing.
See Figure \ref{fig:splicing}.

Starting with any $D$, let $s_A:=s_A(D)$ (resp. $s_B:=s_B(D)$) denote
the crossing--free diagram obtained by applying the $A$--resolution
(resp.  $B$--resolution) to all the crossings of $D$.  We obtain
graphs $\GA$, $\GB$ as follows: The vertices of $\GA$ are in
one-to-one correspondence with the circles of $s_A$.  Every
crossing of $D$ gives rise to two arcs of the $A$--resolution.  These
will each be associated with a component of $s_A$, so correspond to
vertices of $\GA$.  Add an edge to $\GA$ connecting these two vertices
for each crossing of $D$.  In a similar manner, construct the
$B$--graph $\GB$ by considering components of $s_B$.

A link diagram $D$ is called \emph{adequate} if the graphs $\GA$,
$\GB$ contain no edges with both endpoints on the same vertex.  A
link is called \emph{adequate} if it admits an adequate diagram.

Let $v_A$, $e_A$ denote the number of vertices and edges of $\GA$,
respectively.  Similarly, let $v_B$ and $e_B$ denote the number of
vertices and edges of $\GB$.  The reduced graph $\GRA$ is obtained
from $\GA$ by removing multiple edges connected to the same pair of
vertices.  The reduced graph $\GRB$ is obtained similarly.  Let $e'_A$
(resp. $e'_B$) denote the number of edges of $\GRA$ (resp. $\GRB$).  A
proof of the following lemma can be found in \cite{dasbach-lin:volumeish}.

\begin{lemma}[Stoimenow]\label{lem:oxs} 
Let $D$ be an adequate diagram of a link $K$.  Let $\beta$ and $\beta'$
  be the second and next-to-last coefficients of $J_K(t)$.  Then
$$\abs{\beta} + \abs{\beta'} \; = \; e'_A+e'_B- v_A-v_B + 2.$$
\end{lemma}

\subsection{Tangle addition}
Let $D$ be a diagram of a link $K$ obtained by summing strongly
alternating diagrams of tangles $T_1, \dots, T_n$ as in the statement
of Theorem \ref{tangleJP}.  By work of Lickorish and Thistlethwaite
\cite{lith}, $D$ is an adequate diagram; thus the result stated above
applies to $K$.  To estimate the quantity $e'_A+e'_B- v_A-v_B + 2$ we
need to examine the loss of edges as one passes from $\GA$, $\GB$ to
the reduced graphs.

Let $T$ denote a strongly alternating tangle.  Recall $T$ lies inside
a disk on the plane.  One can define the $A$--graph $\Gamma_A(T)$,
and the $B$--graph $\Gamma_B(T)$, corresponding to $T$ in a way
similar to the diagram $D$, by resolving the crossings of $T$ in the
interior of the disk.  Similarly, we can consider the reduced $A$ and
$B$ graphs of $T$; denote them by $\Gamma'_A(T)$ and $\Gamma'_B(T)$,
respectively.

In an alternating diagram of a tangle or link, every component of
$s_A$ and $s_B$ follows along the boundary of a region of the diagram.
Thus the vertices of $\Gamma_A(T)$ and $\Gamma_B(T)$ are in 1--1
correspondence with regions in the diagram of $T$.  These graphs will
have two types of vertices: \emph {interior vertices}, corresponding
to regions that lie entirely in the disk, and two \emph {exterior
vertices}, corresponding to the two regions with sides on the boundary
of the disk.

\begin{lemma}
Let $T$ be an alternating tangle.  Then the only edges lost as we pass
from $\Gamma_A(T)$, $\Gamma_B(T)$ to $\Gamma'_A(T)$, $\Gamma'_B(T)$
are multiple edges from twist regions.  In a twist region with $c_R$
crossings, we lose exactly $c_R-1$ edges.
\label{lemma:tangle-loss}
\end{lemma}

Compare this to \cite{dasbach-futer...,dasbach-lin:volumeish}, where
similar statements are proved for knots and links.

\begin{proof}
We have observed above that the vertices of $\Gamma_A(T)$ and
$\Gamma_B(T)$ are in 1--1 correspondence with regions in the diagram
of $T$.  Thus if edges $e$ and $e'$ connect the same pair of vertices,
the loop $e \cup e'$ passes through exactly two regions of the
diagram, while intersecting the diagram at two crossings.  Therefore,
these crossings are equivalent, and belong to the same twist region.

Conversely, a twist region $R$ with $c_R$ crossings corresponds to a
pair of vertices that are connected by $c_R$ edges.  Therefore, as we
pass to the reduced graphs $\Gamma'_A(T)$ and $\Gamma'_B(T)$, we lose
exactly $c_R-1$ edges from $R$.
\end{proof}

\begin{figure}
\centerline{\includegraphics{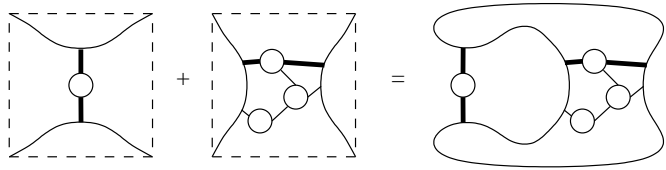}}
\caption{An example of unexpected losses.  Upon tangle addition, the
	dark edges connect the same pair of vertices.}
\label{fig:unexpected}
\end{figure}

As we add several tangles to obtain a link diagram $D$, we may
encounter additional, \emph{unexpected} losses of edges, because the
two exterior vertices in a tangle become amalgamated when we perform
the Conway sum.
Note that because each tangle is chosen to be strongly alternating,
the two exterior vertices of any tangle cannot be connected to each
other by an edge in the tangle.  Thus each edge with an endpoint on
one exterior vertex must have the other endpoint on an interior
vertex.  Then when we do the sum, the only way to pick up an
unexpected loss is to have a tangle with both exterior vertices
connected by edges to the same interior vertex, and then in the sum to
have those two exterior vertices identified to each other.  See Figure
\ref{fig:unexpected}.

\begin{define}
Let $D$ be a diagram obtained by summing strongly alternating diagrams
of tangles $T_1, \dots, T_n$.  Let $\lbigon(D)$ denote the total
loss of edges as we pass from $e_A+e_B$ to $e'_A+e'_B$ which come from
equivalent crossings in the same tangle $T_i$.  Let $\ltop(D)$ denote
the total loss of edges coming from tangle addition.

For a tangle $T\in \{T_1, \ldots, T_n\}$ a \emph{bridge} of
$\Gamma_A(T)$ (resp. $\Gamma_B(T)$) is a subgraph consisting of an
interior vertex $v$, the two exterior vertices $v', v''$ and two edges
$e', e''$ such that $e'$ connects $v$ to $v'$ and $e''$ connects $v$
to $v''$.  The bridge is called \emph{inadmissible} iff $v', v''$
collapse to the same vertex in $\GA$ (resp. $\GB$).
This is the situation of Figure \ref{fig:unexpected}.
\label{def:losses}
\end{define}

It follows that $e_A + e_B - e'_A - e'_B \; = \; \lbigon+ \ltop.$
By Lemma \ref{lemma:tangle-loss}, we have $\lbigon=c(D)-\tw (D)$.  In
the next lemma we estimate $\ltop$.

\begin{lemma}\label{lemma:ltop-estimate}
Let $T_1, T_2$ be strongly alternating tangles
whose Conway sum is a knot diagram $D(K)$.
 Then
$$\ltop(D) \; \leq \; \frac{\twist}{2}+4.$$
\end{lemma}

\begin{proof}
For $T\in \{T_1, T_2 \}$ let $b_A(T)$, $b_B(T)$ denote the number of
bridges in $\Gamma_A(T)$, $\Gamma_B(T)$, respectively.  Then, the
contribution of $T$ to $\ltop$ is at most $b_A(T)$+$b_B(T)$.

Now let $b$ be a bridge of a tangle $T$.  There are two possibilities
for $b$:
\begin{enumerate}
	\item[(I)] The edges $e'$, $e''$ of definition \ref{def:losses} do
		\emph{not} come from resolutions of a single twist region.
	\item[(II)] The edges $e'$, $e''$ of definition \ref{def:losses}
		\emph{do} come from resolutions of a single twist region.
\end{enumerate}

Note for a type (II) bridge, the interior vertex $v$ comes from a
bigon of the diagram, and the corresponding twist region has exactly
two crossings.  This is illustrated in Figure \ref{fig:bridge=2}(a)
for $\Gamma_A(T)$: A type (II) bridge gives two crossings as in that
figure, where shaded regions become vertices of $\Gamma_A(T)$.  By
definition of twist region, there is a simple closed curve meeting the
diagram in exactly the two crossings, as shown by the dotted line.
The strands of the crossing cannot cross the shaded region inside the
dotted line, since this becomes a single vertex of $\Gamma_A(T)$.
Since the diagram is prime, the tangle within the dotted line must be
trivial, consisting of two unknotted arcs.  Finally, no other crossing
can be in the same equivalence class as the two shown, because such a
crossing would have to lie in one of the shaded regions, but these are
vertices of $\Gamma_A(T)$.

\begin{figure}
\begin{center}
	\begin{tabular}{ccccc}
		(a) &
		\includegraphics{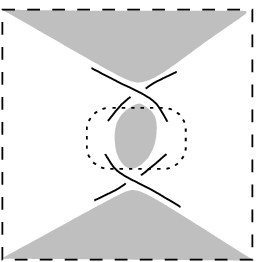}
		\hspace{.5in} &
		(b) &
		\includegraphics{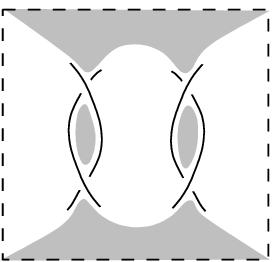}
	\end{tabular}
\end{center}
\caption{(a)  A type (II) bridge gives a bigon of the diagram $D$.
(b) More that one type (II) bridge implies the tangle has more than
	one component.}
\label{fig:bridge=2}
\end{figure}

As we pass from the graphs $\GA$, $\GB$ to the reduced ones $\GRA$,
$\GRB$ each bridge loses exactly one of the edges $e'$, $e''$.  The
contribution to $\ltop$ from type (I) bridges is half of the number of
twist regions in $T$ involved in such bridges.

As for type (II) bridges, if a tangle $T \in \{T_1, T_2\}$ is such
that $\Gamma_A(T)$ or $\Gamma_B(T)$ has more than one bridge of type
(II), then $K$ has more than one component.  This is illustrated in
Figure \ref{fig:bridge=2}(b).  If $\Gamma_A(T)$ has more than one
bridge of type (II), $T$ must be as shown in the figure, with shaded
regions corresponding to vertices of $\Gamma_A(T)$, and possibly
additional crossings in the white regions of the diagram.  Note the
four strands in the center region must connect to form one or two distinct link
components.

Also observe that there cannot simultaneously be two--crossing twist regions
 connecting the east side to the west and the north to the south.
 Hence we may conclude that $T$ contains at most one bridge
of type (II).

\begin{figure}[h]
\centerline{\includegraphics{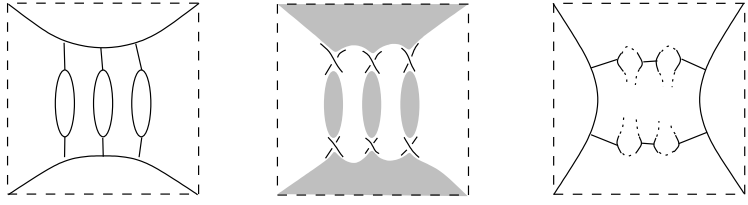}}
\caption{If $\Gamma_A(T)$ contains at least three bridges, as at left,
	then $T$ is as at center, so $\Gamma_B(T)$, at right, cannot contain
	any bridges.}
\label{fig:bridge>2}
\end{figure}

{\it Case 1:} Suppose that $b_A(T)\geq 3$ or $b_B(T)\geq 3$.  Without
loss of generality, say $b_A(T) \geq 3$.  Then we claim that $b_B(T) =
0$.  This is illustrated in Figure \ref{fig:bridge>2}:  If
$\Gamma_A(T)$ contains at least three bridges, then the tangle $T$
must have crossings in the form of the center of that figure.  Note no
edge of $\Gamma_B(T)$ can run through the shaded regions of that
figure, else it will correspond to a crossing which would split an
interior bridge vertex of $\Gamma_A(T)$.  Thus any path from the left
to the right exterior vertex of $\Gamma_B(T)$ must contain at least
three edges, so $\Gamma_B(T)$ cannot contain any bridges.

Then either there are no type (II) bridges in $\Gamma_A(T)$, and then at
most $\tw(T)/2$ bridges, or there are at most $(\tw(T)-1)/2$ bridges
of type (I) and a single bridge of type (II).  In either case, the
contribution of $T$ to $\ltop$ is at most
$$b_A(T)+b_B(T) \:\leq\: \frac{(\tw(T)-1)}{2} +1 
\:<\: \frac{\tw(T)}{2} + 2.$$

\begin{figure}[h]
\centerline{\includegraphics{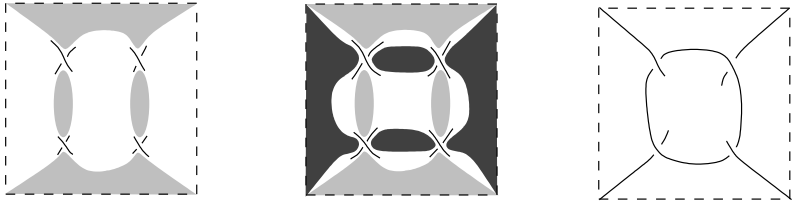}}
\caption{When $b_A(T) = b_B(T) = 2$, $T$ must be as shown in the
center.  This forces $T$ to have at least two components, as at
right.}
\label{fig:bridge4}
\end{figure}

{\it Case 2:} Next, suppose that $b_A(T)\leq 2$ and $b_B(T)\leq 2$.
Then $b_A(T)+b_B(T)\leq 4$, and the maximum contribution to $\ltop(T)$
occurs when $b_A(T)=b_B(T)=2$.  However, we now show that when $b_A(T)
= b_B(T) = 2$, we actually have a link rather than a knot.  This is
illustrated in Figure \ref{fig:bridge4}.  If $b_A(T)=2$, the tangle
diagram must be as in the left of that figure.  Similarly, if
$b_B(T)=2$, the tangle diagram must be as in the left, but rotated 90
degrees.  Since edges of $\Gamma_B(T)$ cannot pass through the
vertices of $\Gamma_A(T)$ (shaded regions of the figure), and vice
versa, the only possibility is that the tangle $T$ has the form in the
center.  Here the lighter shaded regions become vertices of
$\Gamma_A(T)$, and the darker become vertices of $\Gamma_B(T)$.  But
then $T$ must actually have a diagram as on the right of the figure,
because closures of the diagram are prime, implying the diagram is
prime.  Note the tangle must consist of at least two components.

So suppose $b_A(T) = 2$ and $b_A(T) = 1$, or vice versa.  Then because
there is at most one bridge of type (II), $T$ must contain at least
two twist regions.  Thus
$$b_A(T)+b_B(T) = 3 \leq \frac{\tw(T)}{2} + 2.$$

For strongly alternating tangles, the twist number is additive under
tangle addition, which can be seen as follows.  Suppose $T_a$ and
$T_b$ are tangles whose sum has diagram $\Delta$.  First, since
equivalent crossings in a tangle are still equivalent after tangle
addition, and tangle addition does not produce more crossings,
$\tw(\Delta) \leq \tw(T_a) + \tw(T_b)$.  Suppose $\tw(\Delta)$ is
strictly less than $\tw(T_a) + \tw(T_b)$.  That means two twist
regions in distinct tangles become equivalent under tangle sum.  By
definition, there exists a simple closed curve $\gamma$ in $\Delta$
meeting just a crossing in $T_a$, and just a crossing in $T_b$.  It
must run through the unit square bounding $T_a$.  Note by parity,
$\gamma$ either intersects the north and south edges of the unit square,
or the east and west edges.  But in the first case, the denominator
of the tangle is not prime, and in the second the numerator is not
prime, contradicting strongly alternating.  Thus the twist number is
additive under tangle addition.

The previous inequality therefore implies that
$$\ltop(D) \; \leq \; \sum_{i=1,2} (b_A(T_i) + b_B(T_i)) \; \leq \;
\frac{\twist }{2}+4.$$

\vspace{-3ex}
\end{proof}

\begin{proof}[Proof of Theorem  \ref{tangleJP}]

It is well--known that the Jones polynomial of a link remains
invariant under mutation \cite{lickorish:book}.  Thus, for our
purposes, we are free to modify $D$ by mutation.  After mutation we
can assume that the sum of the tangles $T_1+ \ldots + T_n$ is either a
strongly alternating tangle, or it splits in the form $T+T'$ where
each of $T$, $T'$ is strongly alternating and $T+T'$ is not
alternating.  In the former case we have a stronger result: Dasbach
and Lin \cite{dasbach-lin:volumeish} have shown that $\twist =
\abs{\beta} + \abs{\beta'}$.

So now we assume that $D$ is not alternating.  By work of Lickorish
and Thistlethwaite \cite{lith}, $D$ is an adequate diagram; thus the
results stated above apply for $K$.  By Propositions 1 and 5 of  \cite{lith} (see also
\cite{dasbach-futer...}) we have
\begin{equation}
v_A+v_B=c,
\label{eqn:lith}
\end{equation}
where $c:=c(D)$ denotes the crossing number of $D$.  Now, recall that
every edge of $\GA$ or $\GB$ that is lost as we pass to $\GRA$ and
$\GRB$ either comes from a twist region in a tangle, or an edge of an
inadmissible bridge.  The number of edges lost due to twist regions is
$c-t$, where $t= \twist$.  Thus
$$e_A + e_B - e'_A - e'_B \; = \; (c-t) + \ltop.$$

Now by Lemma \ref{lem:oxs}, we have
\begin{eqnarray*}
\abs{\beta} + \abs{\beta'}
&=& e'_A + e'_B - v_A - v_B + 2 \\
&=& (e'_A + e'_B - e_A - e_B) + e_A + (e_B - v_A - v_B) + 2 \\
&=& -(c-t + \ltop) + c +(c- v_A - v_B)+2 \\
&\geq& t - \ltop +2 \hspace{1.71in} \mbox{(by (\ref{eqn:lith}))}\\
&\geq&t- \frac{t}{2}-4+2=\frac{t}{2}-2 \hspace{1in} \mbox{(by Lemma
	\ref{lemma:ltop-estimate})} 
\end{eqnarray*}
The upper bound on $\abs{\beta} + \abs{\beta'}$ was proved in
Proposition 4.6 of \cite{fkp-07}.
\end{proof}

\section*{Acknowledgements}
We thank Peter Milley for writing valuable code to check cusp
area estimates, and for helpful correspondence.  We also thank Bruno
Martelli, Carlo Petronio, and Genevieve Walsh for helpful correspondence.

\bibliographystyle{hamsplain}
\bibliography{biblio.bib}

\end{document}